\documentclass[12pt,psamsfonts]{article}
\usepackage[dvips]{epsfig}
\usepackage{graphics}
\usepackage{amsmath}
\usepackage{amsfonts}
\usepackage{amssymb}
\usepackage{amsthm}
\usepackage{bm}
\usepackage{enumerate}
\usepackage[mathscr]{eucal}
\usepackage{inputenc}
\usepackage[T2A]{fontenc}
\theoremstyle{plain}
\newtheorem{theorem}{Theorem}[section]
\newtheorem{lemma}{Lemma}[section]

\newtheorem{remark}{Remark}[section]

 \numberwithin{equation}{section}

\begin{document}
\title{\textbf{The truncated Fourier operator. II}}
\author{\textbf{Victor Katsnelson} \and \textbf{Ronny Machluf}}
\footnotetext{\hspace*{-4.0ex}\textbf{Mathematics Subject
Classification: (2000).} Primary 47A38; Secondary 47B35, 47B06,
47A10.\endgraf \hspace*{-2.0ex}\textbf{Keywords:} Truncated
Fourier operator, differential operators commuting with the
truncated Fourier operator on finite or semi-infinite interval.}
\date{\ }
 \maketitle
 \abstract{
 For \(E\) being one of the three sets: the whole real axis,
 a finite symmetric interval and the positive semiaxis,
 we discuss the simplest differential operators of the second order
 which commute with the truncated Fourier operator \(\mathscr{F}_E\). }
\setcounter{section}{1}
\section{Operators which commute with the
truncated Fourier operator}

 Let \(F\) be a linear operator acting in a vector space
 \(\mathscr{X}\) over \(\mathbb{C}\), and let \(L\) be another
 linear operator acting in \(\mathscr{X}\). We assume that the
 operator \(L\) commutes with \(F\) that is the equality
 \begin{equation}%
 \label{CoRe}
 FLx=LFx\quad \forall x\in \mathscr{X}
 \end{equation}%
 holds.
 For \(\lambda\in\mathbb{C}\), Let \(\mathcal{X_\lambda}\) be the
 eigenspace of the operator \(L\) corresponding the eigenvalue \(\lambda\):
\begin{equation}%
 \label{EiSp}
\mathscr{X_\lambda}=\lbrace{}x\in\mathscr{X}:\,Lx=\lambda{}x\rbrace\,.
\end{equation}%
From \eqref{CoRe} and \eqref{EiSp} it follows that
\begin{equation}%
 \label{InSu}
 F\mathscr{X_\lambda}\subseteq\mathscr{X_\lambda}\,,
\end{equation}%
that is \(\mathscr{X_\lambda}\) is an invariant subspace of the
operator \(F\).

 Thus the use of operators \(L\) commuting  with the given
operator \(F\) may be a helpful tool for study of the operator
\(F\) if an information on spectral properties of \(L\) is
available.

We apply this reasoning taking the truncated Fourier operator
\(\mathscr{F}_E\) as the operator \(F\).

\hspace*{1.0ex}\textsf{1}.\,\,First of all, we consider the case
when \textsf{the set \(E\) is symmetric}, that is \(E=-E\).

On the symmetric set \(E\), the symmetry transformation
\(t\to{}-t\) acts. This transformation generates the operator
\begin{equation}%
\label{InOp}%
\mathscr{J}x(t)=x(-t)\,,\quad \mathscr{J}:\,L^2(E)\to{}L^2(E)\,.
\end{equation}%
The operator \(\mathscr{J}\) possesses the properties:
\begin{equation}%
\label{PSO}%
\mathscr{J}^2=I,\quad \mathscr{J}=\mathscr{J}^\ast\,,
\end{equation}%
\(I\) in \eqref{PSO}  is the identity operator in \(L^2(E)\). The
spectrum of the operator \(\mathscr{J}\) consists of two points:
\(\lambda=1\) and \(\lambda=-1\), which are  eigenvalues of
\(\mathscr{J}\). We denote the appropriate eigenspaces of
\(\mathscr{J}\) by \(\mathscr{X}_e\) and \(\mathscr{X}_o\):
\begin{subequations}
\label{eo}
\begin{align}
\label{eo1}%
\mathscr{X}_e=\lbrace{}x\in\L^2_E:\,\mathscr{J}x=\phantom{-}x,\ \
&\textup{that
is} \ \ x(t)=\phantom{-}x(-t),\ \ t\in{}E\,;\\
\label{eo2}%
 \mathscr{X}_o=\lbrace{}x\in\L^2_E:\,\mathscr{J}x=-x,\
\ &\textup{that is} \ \ x(t)=-x(-t),\ \ t\in{}E\,;
\end{align}
\end{subequations}
(The subspace \(\mathscr{X}_e\) consists of even functions, the
subspace \(\mathscr{X}_o\)\,---\,of odd.) It is clear that
\begin{equation}
\label{DiS}%
 L^2(E)=\mathscr{X}_e\oplus{}\mathscr{X}_o\,.
\end{equation}
The operator \(\mathscr{F}_E\) commutes with \(\mathscr{J}\):
\begin{subequations}
\label{comrel}%
\begin{equation}
\label{comrel1}%
 \mathscr{J}\mathscr{F}_E=\mathscr{F}_E\mathscr{J}\,,
\end{equation}
and also
\begin{equation}
\label{comrel2}%
 \mathscr{J}\mathscr{F}^{\,\ast}_E=\mathscr{F}^{\,\ast}_E\mathscr{J}\,.
\end{equation}
\end{subequations}
 Therefore,
 \begin{subequations}
 \label{PrPa}
\begin{align}%
\label{PrPa1}
\mathscr{F}_E\mathscr{X}_e\subseteq\mathscr{X}_e,\quad
\mathscr{F}_E\mathscr{X}_o\subseteq\mathscr{X}_o\,;\\
\mathscr{F}^{\ast}_E\mathscr{X}_e\subseteq\mathscr{X}_e,\quad
\mathscr{F}^{\ast}_E\mathscr{X}_o\subseteq\mathscr{X}_o\,.
\end{align}
\end{subequations}
So the pair of complementary subspaces \(\mathscr{X}_e\),
\(\mathscr{X}_o\) reduces each of the operators \(\mathscr{F}_E\)
and \(\mathscr{F}^{\,\ast}_E\):
\begin{subequations}
\label{Red}
\begin{align}
\label{Red1}
\mathscr{F}_E={\mathscr{F}_E}_{|\mathscr{X}_e}\oplus{}{\mathscr{F}_E}_{|\mathscr{X}_o}\,,\\
\label{Red2}
\mathscr{F}^{\ast}_E={\mathscr{F}^{\ast}_E}_{|\mathscr{X}_e}\oplus{}{\mathscr{F}^{\ast}_E}_{|\mathscr{X}_o}\,,
\end{align}
\end{subequations}
where \({\mathscr{F}_E}_{|\mathscr{X}_e}\),
\({\mathscr{F}_E}_{|\mathscr{X}_o}\),
\({\mathscr{F}^{\ast}_E}_{|\mathscr{X}_e}\) and
\({\mathscr{F}^{\ast}_E}_{|\mathscr{X}_o}\) are the restrictions
of the operators \(\mathscr{F}_E\) and \(\mathscr{F}^{\ast}_E\)
onto their invariant subspaces \(\mathscr{X}_e\) and
\(\mathscr{X}_o\) respectively.
 According to Theorem 1.3, the
operator \(\mathscr{F}_E\), as well as the adjoint operator
\(\mathscr{F}^{\ast}_E\), are normal operators in \(L^2(E)\). The
restriction of a normal operator onto its reducing subspace is a
normal operator as well. Therefore the operators
\({\mathscr{F}_E}_{|\mathscr{X}_e}\) and
\({\mathscr{F}_E}_{|\mathscr{X}_o}\) are normal operators acting
in the spaces \(\mathscr{X}_e\) and \(\mathscr{X}_o\)
respectively.

The relation %
\begin{equation}
\label{Jad}%
 \mathscr{F}_E^{\,\ast}=\mathscr{F}_E\mathscr{J}
\end{equation}
is crucial for us. Since
\begin{subequations}
\label{FES}
\begin{align}
\label{FES1} %
\hspace*{6.0ex}
\mathscr{J}_{|\mathscr{X}_e}&=\phantom{-}I_{\mathscr{X}_e},\quad
(I_{\mathscr{X}_e}- \textup{the identity operator in}\ \
\mathscr{X}_e)\,,\\
\label{FES2} %
\hspace*{6.0ex}
\mathscr{J}_{|\mathscr{X}_o}&=-I_{\mathscr{X}_o},\quad
(I_{\mathscr{X}_o}- \textup{the identity operator in}\ \
\mathscr{X}_o)\,,
\end{align}
\end{subequations}
then
\begin{subequations}
\label{SFs}
\begin{align}
\label{SFs1}
(\mathscr{F}_E)_{|\mathscr{X}_e}&=\phantom{-}(\mathscr{F}^{\,\ast}_E)_{|\mathscr{X}_e}\\
\label{SFs2}
(\mathscr{F}_E)_{|\mathscr{X}_o}&=-(\mathscr{F}^{\,\ast}_E)_{|\mathscr{X}_o}\,.
\end{align}
\end{subequations}
It is also clear that
\((\mathscr{F}^{\,\ast}_E)_{|\mathscr{X}_e}=((\mathscr{F}_E)_{|\mathscr{X}_e})^{\,\ast}\,,\quad
(\mathscr{F}^{\,\ast}_E)_{|\mathscr{X}_o}=((\mathscr{F}_E)_{|\mathscr{X}_o})^{\,\ast}\,.\)
(The operator which is adjoint to the restricted operator is equal
to the restriction of the ajoint operator\,.)

 Thus, the operators \((\mathscr{F}_E)_{|\mathscr{X}_e}\) and
\(i^{-1}(\mathscr{F}_E)_{|\mathscr{X}_o}\) are selfadjoint
opearators in the subspaces \(\mathscr{X}_e\) and
\(\mathscr{X}_o\) of even and odd functions from \(L^2(E)\)
respectively. These operators are nothing more then the cosine and
sine Fourier transformation:
\begin{subequations}
\label{pft}
\begin{align}
\label{pft1}
\textup{For}\ x\in\mathscr{X}_e,\ \ F_Ex(t)&=\phantom{i}%
\frac{1}{\sqrt{2\pi}}\int\limits_{E}\cos{(t\xi)}\,x(\xi)\,d\xi\,,\quad
t\in{}E\,,\\
 \label{pft2}
\textup{For}\ x\in\mathscr{X}_o,\ \
F_Ex(t)&=i\frac{1}{\sqrt{2\pi}}\int\limits_{E}\sin{(t\xi)}\,x(\xi)\,d\xi\,,\quad
t\in{}E\,,
\end{align}
\end{subequations}
 Each of the subspaces \(\mathscr{X}_e\) and \(\mathscr{X}_o\)
is invariant with respect to each of operators \(\mathscr{F}^2_E\)
and \(\mathscr{F}^{\ast}_E\mathscr{F}_E\). Taking into account
\eqref{Jad} and \eqref{FES}, we see that
\begin{subequations}
\label{SFf}
\begin{align}
\label{SFf1}%
 ({\mathscr{F}_E}_{|\mathscr{X}_e})^2&=
\phantom{-}(\mathscr{F}^{\,\ast}_E\mathscr{F}_E)_{|\mathscr{X}_e}\\
\label{SFf2}
({\mathscr{F}_E}_{|\mathscr{X}_o})^2&=-(\mathscr{F}^{\,\ast}_E\mathscr{F}_E)_{|\mathscr{X}_o}\,.
\end{align}
\end{subequations}
The operator \(\mathscr{F}^{\,\ast}_E\mathscr{F}_E\), and its
restrictions
\((\mathscr{F}^{\,\ast}_E\mathscr{F}_E)_{|\mathscr{X}_e}\),
\((\mathscr{F}^{\,\ast}_E\mathscr{F}_E)_{|\mathscr{X}_o}\)  onto
invariant subspaces \(\mathscr{X}_e\), \(\mathscr{X}_o\) are
non-negative contractive operators.
 Therefore, the spectra
\(\sigma((\mathscr{F}^{\,\ast}_E\mathscr{F}_E)_{|\mathscr{X}_e})\)
and
\(\sigma((\mathscr{F}^{\,\ast}_E\mathscr{F}_E)_{|\mathscr{X}_o})\)
of the operators
\((\mathscr{F}^{\,\ast}_E\mathscr{F}_E)_{|\mathscr{X}_e}\),
\((\mathscr{F}^{\,\ast}_E\mathscr{F}_E)_{|\mathscr{X}_o}\) are
subsets of the unit interval of the real axis:
\begin{equation}
\sigma((\mathscr{F}^{\,\ast}_E\mathscr{F}_E)_{|\mathscr{X}_e})\subseteq\,[0,\,1]\,,\quad
\sigma((\mathscr{F}^{\,\ast}_E\mathscr{F}_E)_{|\mathscr{X}_o})\subseteq\,[0,\,1]\,.
\end{equation}
According to the spectral mapping theorem (see, for example,
\cite[Chapter I, Theorem 1.3]{GGK}),
\begin{equation}
\label{SMT}
(\sigma(\mathscr{F}_E)_{|\mathscr{X}_e})^2=%
\sigma((\mathscr{F}^{\,\ast}_E\mathscr{F}_E)_{|\mathscr{X}_e}),\quad%
(\sigma(\mathscr{F}_E)_{|\mathscr{X}_o})^2=-%
\sigma((\mathscr{F}^{\,\ast}_E\mathscr{F}_E)_{|\mathscr{X}_o})\,,
\end{equation}
hence
\begin{equation}
\label{SMTp}%
 (\sigma(\mathscr{F}_E)_{|\mathscr{X}_e})^2\subseteq%
\sigma((\mathscr{F}^{\,\ast}_E\mathscr{F}_E)),\quad
-(\sigma(\mathscr{F}_E)_{|\mathscr{X}_o})^2\subseteq%
\sigma((\mathscr{F}^{\,\ast}_E\mathscr{F}_E))\,.
\end{equation}
In particular,
\begin{equation}%
\label{SMTpp}%
\sigma(({\mathscr{F}_E})_{|\mathscr{X}_e})\subseteq[-1,\,1],\quad%
\sigma(({\mathscr{F}_E})_{|\mathscr{X}_o})\subseteq[-i,\,i]
\end{equation}
(In \eqref{SMTpp},
\(\,[-i,\,i]=\lbrace{}z=iy\in\mathbb{C}:\,-1\leq{}y\leq{}1\rbrace\).)

So, we proof the following
\begin{theorem}
\label{CSS} Assume that the set \(E\) is symmetric.

\noindent%
Then: \\
\hspace*{2.0ex}\textup{1.}
\begin{minipage}[t]{0.92\linewidth}
The space \(L^2(E)\) is the orthogonal sum of the subspaces
\(\mathscr{X}_e\) and \(\mathscr{X}_o\) of even and odd functions,
\eqref{eo}: \( L^2(E)=\mathscr{X}_e\oplus{}\mathscr{X}_o\,.\)
\end{minipage}\\
\hspace*{2.0ex}\textup{2.}
\begin{minipage}[t]{0.92\linewidth}
Each of the subspaces \(\mathscr{X}_e\) and \(\mathscr{X}_o\) is
invariant with respect the operators \(\mathscr{F}_E\)\,,
\[\mathscr{F}_E\mathscr{X}_e\subseteq\mathscr{X}_e,\quad
\mathscr{F}_E\mathscr{X}_o\subseteq\mathscr{X}_o\,.\] The pair of
complementary subspaces \(\mathscr{X}_e\) and \(\mathscr{X}_o\)
reduces the operators \(\mathscr{F}_E\):
\[\mathscr{F}_E={\mathscr{F}_E}_{|\mathscr{X}_e}\oplus{}{\mathscr{F}_E}_{|\mathscr{X}_o}\,,\]
where \({\mathscr{F}_E}_{|\mathscr{X}_e}\) are the restrictions of
the operator \(\mathscr{F}_E\) onto its invariant subspaces
\(\mathscr{X}_e\) and \(\mathscr{X}_o\) respectively.

The operators \({\mathscr{F}_E}_{|\mathscr{X}_e}\) and
\({\mathscr{F}_E}_{|\mathscr{X}_o}\) are integral transforms:
cosine and sine Fourier transforms respectively, \eqref{pft}.
\end{minipage}\\
\hspace*{2.0ex}\textup{3.}
\begin{minipage}[t]{0.93\linewidth}
The operators \({\mathscr{F}_E}_{|\mathscr{X}_e}\) and
\(i^{-1}{\mathscr{F}_E}_{|\mathscr{X}_o}\) are selfadjoint
contractive operators acting in the spaces \(\mathscr{X}_e\) and
\(\mathscr{X}_e\) respectively.
\end{minipage}\\
\hspace*{2.0ex}\textup{4.}
\begin{minipage}[t]{0.93\linewidth}
The operator \(\mathscr{F}_E\) is a normal operator, and the
operators \({\mathscr{F}_E}_{|\mathscr{X}_e}\oplus{}0\) and
\(0\oplus{\mathscr{F}_E}_{|\mathscr{X}_o}\) are real and imaginary
parts of \(\mathscr{F}_E\).
\end{minipage}\\
\hspace*{2.0ex}\textup{5.}
\begin{minipage}[t]{0.93\linewidth}
The spectrum \(\sigma(\mathscr{F}_E)\) of the operator
\(\mathscr{F}_E\) is a subset of the cross \(Q\), where
\(Q=[-1,1]\cup[-i,i].\)\\ More precisely, \eqref{SFf}, if point
\(\lambda\in{}Q\) belongs to the spectrum
\(\sigma(\mathscr{F}_E)\) of the operator \(\mathscr{F}_E\), then
the point \(|\lambda|^2\) belongs to the spectrum
\(\sigma(\mathscr{F}^{\ast}_E\mathscr{F}_E)\) of the operator
\(\mathscr{F}^{\ast}_E\mathscr{F}_E)\). Conversely, if the point
\(\rho^2\), where \(0\leq\rho\leq{}1\) , belongs to the spectrum
\(\sigma(\mathscr{F}^{\ast}_E\mathscr{F}_E)\) of the operator
\(\mathscr{F}^{\ast}_E\mathscr{F}_E)\), then at least one of the
four points \(i^{k}\rho,\,k=0,\,1,\,2,\,3,\) belongs to the
spectrum \(\sigma(\mathscr{F}_E)\) of the operator
\(\mathscr{F}_E\).
\end{minipage}\\
\hspace*{2.0ex}\textup{6.}
\begin{minipage}[t]{0.93\linewidth}
If the point \(\rho^2\) is an eigenvalue of the operator
\(\mathscr{F}^{\,\ast}_E\mathscr{F}_E\) of multiplicity
\(\kappa(\rho^2,\,\mathscr{F}^{\,\ast}_E\mathscr{F}_E)\), then at
least one of the four points \(i^{k}\rho,\,k=0,\,1,\,2,\,3,\) is
an eigenvalue of the
operator \(\mathscr{F}_E\), and the total multiplicity %
\(\sum\limits_{0\leq{}k\leq{}3}\kappa(i^{k}\rho,\,\mathscr{F}_E)\)
of these eigenvalues is equal to
\(\kappa(\rho^2,\,\mathscr{F}^{\,\ast}_E\mathscr{F}_E)\).\\
In particular, if the point \(\rho^2\) is a simple eigenvalue of
the operator \(\mathscr{F}^{\,\ast}_E\mathscr{F}_E\), then
precisely one of the four points \(i^{k}\rho,\,k=0,\,1,\,2,\,3,\)
is an eigenvalue of the operator \(\mathscr{F}_E\) and this
eigenvalue is simple.
\end{minipage}
\end{theorem}

In fact, our Theorem \ref{CSS} puts in order what is presented in
\cite[Section III]{Sl1}.
 As we shall see below, if the set \(E\) is
not symmetric, the spectral properties of the operator
\(\mathscr{F}_E\) may be quite different.

\hspace*{1.0ex}\textsf{2.} Less trivial examples of operators
which commute with the truncated Fourier operator are differential
operators. We restrict ourself to consideration of
\emph{selfadjoint} differential operators only.
 The spectrum of an ordinary differential operator is of
finite multiplicity, which does not exceed the order of the
operator. If a selfadjoint differential operator \(\mathscr{L}\)
commutes with the truncated Fourier operator \(\mathscr{F}_E\),
then every eigenspace \(\mathscr{X}_\lambda\) of the operator
\(\mathscr{L}\) is invariant with respect to \(\mathscr{F}_E\). If
moreover the spectrum of \(\mathscr{L}\) is discrete (i.e. if the
system of its eigenvector is complete), then we find the complete
orthogonal system
\(\lbrace\mathscr{X}_\lambda\rbrace_{\lambda\in\sigma(\mathscr{L})}\)
of finite dimensional subspaces each of them is invariant with
respect to \(\mathscr{F}_E\). (And the dimension of each
\(\mathscr{X}_\lambda\) does not exceed the order of the
differential operator \(\mathscr{L}\)). The further spectral
analysis of the operator \(\mathscr{F}_E\) is reduced to the
spectral analysis of the restriction
\((\mathscr{F}_E)_{|\mathscr{X}_\lambda}\) of \(\mathscr{F}_E\)
onto \(\mathscr{L}\). If moreover the spectrum of \(\mathscr{L}\)
is simple, i.e. if \(\dim{}\mathscr{X}_{\lambda}=1\) for every
\(\lambda\) (\(\lambda\) are pairwise different), then the system
of the appropriate eigenvectors of \(\mathscr{L}\) forms also the
basis of eigenvectors of  \(\mathscr{F}_E\).

Usually the study of properties of  eigenfunctions (such as its
asymptotic behavior) of differential operators may be done easier
and in more detail than the study of this properties for
eigenfunctions of integral operator if this integral operator has
no singularities. Therefore the availability of a differential
operators \(\mathscr{L}\) which commutes with \(\mathscr{F}_E\)
may be very useful for the spectral analysis of \(\mathscr{F}_E\).
However, the existence of such \(\mathscr{L}\) is a \emph{lucky
accident}, and to our best knowledge there is no regular way to
search for such \(\mathscr{L}\).

Until the present time, the differential operator \(\mathscr{L}\)
were discussed mainly  which commutes with the operator
\(\mathscr{F}_E^{\ast}\mathscr{F}_E\) rather with the operator
\(\mathscr{F}_E\), that is
\[
\mathscr{L}(\mathscr{F}_E^{\ast}\mathscr{F}_E)=
(\mathscr{F}_E^{\ast}\mathscr{F}_E)\mathscr{L}\,.
\] %
For the set \(E=(-\infty,\infty)\), the differential operator
generated by the formal differential expression
\begin{equation}%
\label{wa}%
(L_Ex)(t)=-\dfrac{d^2x(t)\,}{dt^2}+t^2x(t)
\end{equation}%
commutes with \(\mathscr{F}_E\). This fact has been known and has
been exploited since long ago.

 For \(E\) being a finite
symmetric interval, say \(E=(-a,a),\,\,0<a<\infty\), the
differential operator which commutes with \(\mathscr{F}_E\), has
been discovered in \cite{SlPo} and has been exploited in
\cite{LaP1}, \cite{LaP2}, and since in further works. This
operator is generated by the differential expression
\begin{equation}%
\label{fi}%
(L_Ex)(t)=-\dfrac{d\,\,}{dt}\bigg(1-
\frac{t^2}{a^2}\bigg)\dfrac{dx(t)}{dt}+t^2x(t)\,.
\end{equation}%
and by certain boundary conditions which will be described later.
It is not inconceivable that the fact that the differential
operator generated by the the differential expression \eqref{fi}
commutes with the operator
\(\mathscr{F}_E^{\ast}\mathscr{F}_E,\,\,E=[-a,a]\), was known
before the the serie of work \cite{SlPo}). But surely it is the
work \cite{SlPo} where this fact was emphasized and was firstly
used for the spectral analysis of the operator
\(\mathscr{F}_E^{\ast}\mathscr{F}_E\). %
(In fact, this differential operator commutes not only with the
operator \(\mathscr{F}_E^{\ast}\mathscr{F}_E\), but also with the
operator \(\mathscr{F}_E,\,\,E=[-a,a]\).) In \cite{Sl1} the
multidimensional generalizations of the results obtained in
\cite{SlPo}, \cite{LaP1}, \cite{LaP2} were obtained. In
\cite{Sl2}, the discrete analogs of some results of \cite{SlPo},
\cite{LaP1}, \cite{LaP2} were considered. However, the main break
through was put into effect in \cite{SlPo}, \cite{LaP1},
\cite{LaP2}. For sets \(E\) which are different from the whole
real axes or from finite symmetric interval no differential
operators commuting with \(\mathscr{F}_E\) or
\(\mathscr{F}_E^\ast\mathscr{F}_E\) are found. In \cite{Wal}, the
related problem was considered concerning differential which
commutes with the indicator function of the set \(E\).

There is lot of serendipity in the fact that such differential
operator exists and was found in the case
\(E=[-a,a],\,a\in(0,\infty)\). The common feeling is  that\,%
\footnote{Here we quote from \cite[page 379]{Sl2}} %
\emph{"there is something deeper here that we currently
understand\,--\,%
that there is a way of viewing these problems more abstractly that
will explain heir elegant solution in a more natural and profound
way, so that these nice results will not appear so much as lucky
accident."} Starting from early eighties, F.A.\,Gr\"unbaum takes
many efforts to explain the phenomenon from a more general point
of view. In particular, he related the phenomenon of the existence
of a differential operator commuting with certain integral
operators to the phenomenon of bispectrality. The bispectrality
phenomenon lies in the fact that certain differential (or
difference) operators admit the family of eigenfunctions
satisfying some differential (or difference) equation with respect
to the spectral parameter. See \cite{DuGr}, \cite{Gru1},
\cite{Gru2}, and other works of F.A.\,Gr\"unbaum. There is a
literature devoted to the bispectrality, see for example
\cite{HaKa}.

However, these efforts, being interesting in its own right, did
not furnish the desired result up to now.

For \(E=[0,+\infty)\), the differential operator generated by the
differential expression
\begin{equation}%
\label{ha}%
(L_Ex)(t)=-\dfrac{d\,\,}{dt}\bigg(t^2\,\frac{dx}{dt}\bigg)
\end{equation}%
commutes with both operator \(\mathscr{F}_E\) and
\(\mathscr{F}_E^{\,\,\ast}\).

\begin{remark}
The "factor" \(\Big(1-\frac{t^2}{a^2}\Big)\) , as well as the
"factor" \(t^2\), appearing in the expression \eqref{fi} and
\eqref{ha} respectively,
 play the following role. Checking
the commutativity of the differential operator generated by the
differential expression \(L_E\), \eqref{fi} and \eqref{ha}
respectively, we have to integrate by parts. The above mentioned
factors ensure that the terms outside the integral vanish.
Moreover, in order to these terms vanish the appropriate boundary
conditions on functions from the domain of definition of the
differential operator generated by \(L\) should be imposed.
\end{remark}

Let us check the commutativity condition for the differential
expressions \(L_E\), \eqref{wa}, \eqref{fi}, \eqref{ha}, related
to the Fourier operators \(\mathscr{F}_E\) truncated on the sets
\(E=(-\infty,\infty),\,E=[-a,a]\), and \(E=[0,+\infty)\)
respectively. (In fact, in the first case, the operator
\(\mathscr{F}_E\) is `non-truncated' Fourier operator
\(\mathscr{F}\).) In the present section, we are not concerned
with the detailed description of the domains, where the
appropriate operators are defined. Here we only carry out the
formal integration by parts and compute the terms outside the
integral in every of this three cases of \(E\):
\(E=(-\infty,\infty),\,E=[-a,a]\), and \(E=[0,+\infty)\). The
interpretation of the vanishing condition for the terms outside
the integral as  boundary condition imposed on functions from the
domain of definition of the appropriate differential operator will
be discussed in the next section.

Let us calculate the difference
\(\mathscr{F}_EL_Ex-L_E\mathscr{F}_Ex\) in the above mentioned
three cases: \(E=(-\infty,\infty),\,E=[-a,a]\), and
\(E=[0,+\infty)\)\\[1.0ex]
\textsf{I.\ \ The case \(\boldsymbol{E=(-\infty,\infty)}\).}\\
Let us transform the expression \(\mathscr{F}_EL_Ex\):
\begin{equation*}
(\mathscr{F}_EL_Ex)(t)=\int\limits_{-\infty}^{\infty}%
\Bigg(-\frac{d^2x}{d\xi^2}+\xi^2x(\xi)\Bigg)\,e^{it\xi}d\xi
\end{equation*}
Integrating  twice by parts, we obtain
\begin{equation}
\label{fip}
\int\limits_{a}^{b}\Bigg(-\frac{d^2x(\xi)}{d\xi^2}\Bigg)\,e^{it\xi}d\xi=
-\frac{dx(\xi)}{d\xi}\,e^{it\xi}%
\bigg|_{\xi=a}^{\xi=b}+it\,x(\xi)\,e^{it\xi}\bigg|_{\xi=a}^{\xi=b}+
t^2\int\limits_{a}^{b}x(\xi)\,e^{it\xi}\,d\xi\,.
\end{equation}
The integration by parts is possible under the condition:\\[2.0ex]
\hspace*{2.0ex}\begin{minipage}[h]{0.9\linewidth} {\emph{ The
function \(x(\xi)\) is differentible on every finite interval
\mbox{\((a,b)\subset(-\infty,\infty)\),} and its derivative
\(\frac{dx(\xi)}{d\xi}\) is absolutely continuous on every finite
interval.}}
\end{minipage}\\[2.0ex]
Now, assuming that
\begin{equation}
\label{sc1} %
\int\limits_{-\infty}^{\infty}|x(\xi)|\,d\xi<\infty,\quad
\int\limits_{-\infty}^{\infty}\bigg|\frac{d^2x(\xi)}{d\xi^2}\bigg|\,d\xi<\infty
\end{equation}
 and moreover that
\begin{equation}
\label{bc}%
\lim_{a\to-\infty}x(a)=0,\,\lim_{a\to-\infty}\frac{dx(\xi)}{d\xi}\Big|_{\xi=a}=0,\,
\lim_{b\to\infty}x(b)=0,\,\lim_{b\to
\infty}\frac{dx(\xi)}{d\xi}\Big|_{\xi=b}=0,
\end{equation}
we can pass to the limit in \eqref{fip} as
\(a\to-\infty,\,b\to\infty\). We obtain%
\begin{subequations}
\label{coc}
\begin{equation}
\label{coc1}
\int\limits_{-\infty}^{\infty}\Bigg(-\frac{d^2x(\xi)}{d\xi^2}\Bigg)\,e^{it\xi}d\xi=
t^2\int\limits_{-\infty}^{\infty}x(\xi)\,e^{it\xi}\,d\xi\,\ \
\text{for all}\ \ t\in(-\infty,\infty)\ \ .
\end{equation}
From the other hand,
\begin{equation}
\label{coc2}
\int\limits_{-\infty}^{\infty}\xi^2x(\xi)\,e^{it\xi}\,d\xi=-\frac{d^2\,}{dt^2}
\int\limits_{-\infty}^{\infty}x(\xi)\,e^{it\xi}\,\ \ \ \ \text{for
all}\ \ t\in(-\infty,\infty)
\end{equation}
\end{subequations}
assuming that
\begin{equation*}
\label{scp2}
\int\limits_{-\infty}^{\infty}|x(\xi)|\,d\xi<\infty,\quad
\int\limits_{-\infty}^{\infty}\xi^2|x(\xi)|\,d\xi<\infty\,
\end{equation*}
\begin{lemma}
\label{cwal}
 Let \(x(\xi)\) be a function defined on \((-\infty,\infty)\).
\end{lemma} Assume that  the conditions
\begin{equation}%
\label{msii}%
\textup{a).}\  \displaystyle
\int\limits_{-\infty}^{\infty}\xi^2|x(\xi)|\,d\xi<\infty\,\hspace*{3.0ex}\textup{b).}\
\int\limits_{-\infty}^{\infty}\bigg|\frac{d^2x(\xi)}{d\xi^2}\bigg|\,d\xi<\infty
\end{equation}
hold. Then
\begin{equation}%
\label{zmc}%
 \int\limits_{-\infty}^{\infty}|x(\xi)|\,d\xi<\infty,
\end{equation}
and the condition \eqref{bc} holds.
\begin{proof} To prove \eqref{zmc}, it is enough to prove that
\(\int\limits_{-1}^{1}|x(\xi)|d\xi<\infty\). Since for every
\(\alpha,\,\beta:\,-\infty<\alpha,\beta<\infty\), the equality
\(\displaystyle\frac{dx}{d\xi}(\beta)-\frac{dx}{d\xi}(\alpha)=
\int\limits_{\alpha}{\beta}\frac{d^2x(\xi)}{d\xi}d\xi\) holds, it
follows from the the conditions (\ref{msii}.b) that there exist
the limits
\(\displaystyle\lim_{a\to\-\infty}\dfrac{dx(\xi)}{d\xi}\Big|_{\xi=a}\),
\(\lim_{b\to\infty}\dfrac{dx(\xi)}{d\xi}\Big|_{\xi=b}\). From the
 the condition (\ref{msii}.a) it follows that these limits are
equal to zero. Therefore,
\[\frac{dx(t)}{dt}=\int\limits_{-\infty}^{t}\frac{d^2x(\xi)}{d\xi}d\xi=
-\int\limits_{t}^{\infty}\frac{d^2x(\xi)}{d\xi}d\xi\,\] In
particular,
\[\left|\frac{dx(t)}{dt}\right|\leq{}M,\ \,\textup{where}\ \
M=\int\limits_{-\infty}^{\infty}\bigg|\frac{d^2x(\xi)}{d\xi^2}\bigg|\,d\xi\,,\]
and the function \(x(t)\) satisfy the Lipschitz condition.
\begin{equation*}%
%\label{LC}%
 |x(t_2)-x(t_1)|\leq{}M|t_2-t_1|\quad
\forall\,t_1,\,t_2\,.
\end{equation*}
All the more, \(\int\limits_{-1}^{1}|x(\xi)|d\xi<\infty\). From
(\ref{msii}.a) and the Lipschitz condition it follows that
\(\lim_{t\to\pm\infty}x(t)=0\)\,.
\end{proof}
The above reasoning proves the following
\begin{theorem}
\label{cwa}%
 Let \(x(\xi)\) be a function defined on \((-\infty,\infty)\).
 Assume that the following conditions are satisfied:\\
\hspace*{2.0ex}\textup{1.}
\begin{minipage}[t]{0.92\linewidth}
The function \(x(\xi)\) is differentible on \((-\infty,\infty)\)
and its derivative \(\frac{dx(\xi)}{d\xi}\) is absolutely
continuous on every finite interval
\((a,b)\subset(-\infty,\infty)\).
\end{minipage}\\[1.0ex]
\hspace*{2.0ex}\textup{2.}
\begin{minipage}[t]{0.92\linewidth}
\centerline{\(\displaystyle
\int\limits_{-\infty}^{\infty}\xi^2|x(\xi)|\,d\xi<\infty\,\hspace*{3.0ex}
\int\limits_{-\infty}^{\infty}\bigg|\frac{d^2x(\xi)}{d\xi^2}\bigg|\,d\xi<\infty\)}
\end{minipage}\\[1.0ex]

Then
\(\displaystyle\int\limits_{-\infty}^{\infty}|x(\xi)|\,d\xi<\infty\),
the function \(y(t)\),
\begin{equation}%
\label{FTii}
y(t)=\int\limits_{-\infty}^{\infty}x(\xi)\,e^{it\xi}\,d\xi,
\end{equation}
is twice continuously  differentiable on \((-\infty,\infty)\), and
\begin{equation}
\label{CRii} -\frac{d^2y(t)}{dt^2}+t^2y(t)=
\int\limits_{-\infty}^{\infty}
\bigg(-\frac{d^2x(\xi)}{d\xi^2}+\xi^2x(\xi)\bigg)\,e^{it\xi}d\xi\,\quad
\forall\,t\in(-\infty,\infty).
\end{equation}
\end{theorem}\ {} \\[1.0ex]
\textsf{II.\ \ The case \(\boldsymbol{E=(-a,a),\,0<a<\infty}\).}\\
\hspace*{2.0ex}Let us transform the expression
\(\mathscr{F}_EL_Ex\):
\begin{equation*}
(\mathscr{F}_EL_Ex)(t)=\int\limits_{-a}^{a}%
\left(-\frac{d\,\,}{d\xi}\Bigg(\bigg(1-\frac{\xi^2}{a^2}\bigg)
\frac{dx(\xi)}{d\xi}\Bigg) +\xi^2x(\xi)\right)e^{it\xi}d\xi\,.
\end{equation*}
Integrating  by parts twice, we obtain
\begin{multline}
\label{ffp}%
 \int\limits_{-a+\varepsilon}^{a-\varepsilon}
\left(-\frac{d\,\,}{d\xi}\Bigg(\bigg(1-\frac{\xi^2}{a^2}\bigg)\right)
\frac{dx(\xi)}{d\xi}\Bigg)e^{it\xi}d\xi=\\[1.0ex]
=-\bigg(1-\frac{\xi^2}{a^2}\bigg)\frac{dx(\xi)}{d\xi}\,e^{it\xi}%
\bigg|_{\xi=-a+\varepsilon}^{\xi=a-\varepsilon}
+it\bigg(1-\frac{\xi^2}{a^2}\bigg)x(\xi)%
e^{it\xi}\bigg|_{\xi=-a+\varepsilon}^{\xi=a-\varepsilon}-\\[1.0ex]
-it
\int\limits_{-a+\varepsilon}^{a-\varepsilon}x(\xi)\,\frac{d\,\,}{d\xi}
\left(\bigg(1-\frac{\xi^2}{a^2}\bigg)e^{it\xi}\right)\,d\xi\,.
\end{multline}
Let us assume that
\begin{equation}
\label{sc2} %
\int\limits_{-a}^{a}\bigg|%
\frac{d\,}{d\xi}\bigg(\bigg(1-\frac{\xi^2}{a^2}\bigg)\frac{dx(\xi)}{d\xi}\bigg)\bigg|\,d\xi<\infty
\end{equation}
 and moreover that
\begin{equation}
\label{scb2}
\lim_{\xi\to{}-a+0}(\xi+a)\,\frac{dx(\xi)}{d\xi}=0,\quad
\lim_{\xi\to{}a-0}(\xi-a)\,\frac{dx(\xi)}{d\xi}=0.
\end{equation}

We shall show, (see Lemma \ref{cfin} below), that under the
condition \eqref{sc2}, the conditions
\begin{equation}
\label{Sumfi}%
 \int\limits_{-a}^{a}\big|x(\xi)|\,d\xi<\infty\,
\end{equation}
and
\begin{equation}
\label{scb3}%
 \lim_{\xi\to{}-a}\,(\xi+a)\cdot{}x(\xi)=0,\quad
\lim_{\xi\to{}a}\,(\xi-a)\cdot{}x(\xi)=0.
\end{equation}
hold. Hence, we can pass to the limit as \(\varepsilon\to+0\) in
\eqref{ffp}:
\begin{multline}
\label{ptl}%
 \int\limits_{-a}^{a}
\left(-\frac{d\,\,}{d\xi}\Bigg(\bigg(1-\frac{\xi^2}{a^2}\bigg)\right)
\frac{dx(\xi)}{d\xi}\Bigg)e^{it\xi}d\xi=\\[1.0ex]
= -it \int\limits_{-a}^{a}x(\xi)\,\frac{d\,\,}{d\xi}
\left(\bigg(1-\frac{\xi^2}{a^2}\bigg)e^{it\xi}\right)\,d\xi\,.
\end{multline}
Transforming the integral in the right hand side of \eqref{ptl}
and denoting
\begin{math}
 y(t)=\int\limits_{-a}^{a}x(\xi)e^{it\xi}\,d(\xi),
\end{math}
 we obtain
\begin{multline*}%
-it \int\limits_{-a}^{a}x(\xi)\,\frac{d\,\,}{d\xi}
\left(\bigg(1-\frac{\xi^2}{a^2}\bigg)e^{it\xi}\right)\,d\xi=
t^2y(t)-\frac{d\,\,}{dt}\Bigg(\bigg(1-\frac{t^2}{a^2}\bigg)\frac{dy(t)}{dt}+
\\[1.0ex]
+\frac{d^2y(t)}{dt^2}=
-\frac{d\,\,}{dt}\Bigg(\bigg(1-\frac{t^2}{a^2}\bigg)\frac{dy(t)}{dt}+t^2y(t)-
\int\limits_{-a}^{a}\xi^2x(\xi)e^{it\xi}\,d\xi\,.
\end{multline*}
\begin{lemma}%
\label{cfin}%
Let \(x(\xi)\) be a function defined on the open interval
\((-a,a)\). We assume that:
\\
\hspace*{2.0ex}\textup{1.}
\begin{minipage}[t]{0.92\linewidth}
The function \(x(\xi)\) is differentiable on \((-a,a)\) and that
its derivative \(\dfrac{dx(\xi)}{d\xi}\) is absolutely continuous
on every compact subinterval of the interval \((-a,a)\).
\end{minipage}\\[1.0ex]
\hspace*{2.0ex}\textup{2.}
\begin{minipage}[t]{0.92\linewidth}
The condition \eqref{sc2} is satisfied.
\end{minipage}\\[0.1ex]

Then, the condition
\begin{equation}%
\label{LoB}%
 |x(\xi)|=O\Big(\big|\ln|\xi+a|\big|\Big)\ \ \text{as}
\ \xi\to-a,\ \ \ |x(\xi)|=O\Big(\big|\ln|\xi-a|\big|\Big)\ \
\text{as} \ \xi\to{}a,
\end{equation}
holds, the function \(x(\xi)\) is summable  on the interval
\((-a,a)\), \eqref{Sumfi}, and satisfy the boundary condition
\eqref{scb3}.
\end{lemma}%
\begin{proof} Let
\[-a<\alpha<\beta<a\,.\]
Since
\[\bigg(1-\frac{\xi^2}{a^2}\bigg)\frac{dx(\xi)}{d\xi}\bigg|_{\xi=\beta}\!\!\!-
\bigg(1-\frac{\xi^2}{a^2}\bigg)\frac{dx(\xi)}{d\xi}\bigg|_{\xi=\alpha}\!\!\!=
\int\limits_{\alpha}^{\beta}%
\frac{d\,}{d\xi}\bigg(\bigg(1-\frac{\xi^2}{a^2}\bigg)
\frac{dx(\xi)}{d\xi}\bigg)\,d\xi,\] it follows from \eqref{sc2}
that there exist finite limits
\[\lim_{\alpha\to-a+0}%
\bigg(1-\frac{\xi^2}{a^2}\bigg)\frac{dx(\xi)}{d\xi}\bigg|_{\xi=\alpha},\quad%
\lim_{\beta\to{}a-0}%
\bigg(1-\frac{\xi^2}{a^2}\bigg)\frac{dx(\xi)}{d\xi}\bigg|_{\xi=\beta}\]
In particular,%
\[\frac{dx(\xi)}{d\xi}=O\bigg(\frac{1}{|\xi+a|}\bigg)
\ \text{as} \ \xi\to-a,\ \
\frac{dx(\xi)}{d\xi}=O\bigg(\frac{1}{|\xi+a|}\bigg) \ \text{as} \
\xi\to{}a.
\]
Integrating, be obtain \eqref{LoB}. %
The conditions \eqref{Sumfi} \eqref{scb3} are direct consequences
of \eqref{LoB}.
\end{proof}
 The above reasoning prove the following
\begin{theorem}
\label{cfi}%
Let \(x(\xi)\) be a function defined on \((-a,a)\).
 Assume that the following conditions are satisfied:\\
 \hspace*{2.0ex}\textup{1.}
\begin{minipage}[t]{0.92\linewidth}
The function \(x(\xi)\) is differentible on \((-a,a)\) and its
derivative \(\dfrac{dx(\xi)}{d\xi}\) is absolutely continuous on
every compact subinterval of the interval \((-a,a)\).
\end{minipage}\\[1.0ex]
\hspace*{2.0ex}\textup{2.}
\begin{minipage}[t]{0.92\linewidth}
\centerline{
\(\displaystyle\int\limits_{-a}^{a}\bigg|%
\frac{d\,}{d\xi}\bigg(\bigg(1-\frac{\xi^2}{a^2}\bigg)%
\frac{dx(\xi)}{d\xi}\bigg)\bigg|\,d\xi<\infty\).}
\end{minipage}\\[1.0ex]
\hspace*{2.0ex}\textup{3.}
\begin{minipage}[t]{0.92\linewidth}
The boundary conditions
\begin{equation}
\label{bccfi} %
\displaystyle
\lim_{\xi\to{}-a+0}(\xi+a)\,\frac{dx(\xi)}{d\xi}=0,\quad
\lim_{\xi\to{}a-0}(\xi-a)\,\frac{dx(\xi)}{d\xi}=0\,.
\end{equation}
hold.
\end{minipage}\\[1.0ex]
Then
\begin{equation}
\label{sumc}%
 \int\limits_{-a}^{a}\big|x(\xi)|\,d\xi<\infty\,,
\end{equation}
the function \(y(t)\),
\begin{equation}%
\label{FTiib}%
 y(t)=\int\limits_{-a}^{a}x(\xi)\,e^{it\xi}\,d\xi,
\end{equation}
is defined for every \(t\in\mathbb{C}\) and is an entire function
of exponential type. The following commutational relation holds
\begin{multline}
-\frac{d\,\,}{dt}\Bigg(\bigg(1-\frac{t^2}{a^2}\bigg)\frac{dy(t)}{dt}\Bigg)\,+\,t^2y(t)=\\[1.0ex]
 =\int\limits_{-a}^{a}
\left(-\frac{d\,\,}{d\xi}\bigg(\bigg(1-\frac{\xi^2}{a^2}\bigg)
\frac{dx(\xi)}{d\xi}\bigg)\,+\,\xi^2x(\xi)\right)e^{it\xi}d\xi\,,\quad
\forall\,t\in\mathbb{C}\,.
\end{multline}
\end{theorem}\ {} \\[1.0ex]
\textsf{III.\ \ The case \(\boldsymbol{E=(0,\infty)}\,.\)}\\
\hspace*{2.0ex}Let us transform the expression
\(\mathscr{F}_EL_Ex\):
\begin{equation*}
(\mathscr{F}_EL_Ex)(t)=\int\limits_{0}^{\infty}%
\left(-\frac{d\,\,}{d\xi}\bigg(\xi^2 \frac{dx(\xi)}{d\xi}\bigg)
\right)\,e^{it\xi}d\xi\,.
\end{equation*}
Integrating  twice by parts, we obtain
\begin{multline}
\label{siip}
\int\limits_{\varepsilon}^{N}\left(-\frac{d\,\,}{d\xi}\bigg(\xi^2
\frac{dx(\xi)}{d\xi}\bigg)\right)e^{it\xi}d\xi=\\[1.0ex]
-\xi^2\frac{dx(\xi)}{d\xi}\,e^{it\xi}%
\bigg|_{\xi=a}^{\xi=b}+it\,\xi^2x(\xi)\,e^{it\xi}\bigg|_{\xi=\varepsilon}^{\xi=N}-
it\int\limits_{\varepsilon}^{N}x(\xi)\Bigg(\frac{d\,}{d\xi}\,\bigg(\xi^2e^{it\xi}\bigg)\Bigg)\,d\xi\,.
\end{multline}
Let us assume that
\begin{equation}
\label{ssiic2} %
\textup{a).}\ %
\int\limits_{0}^{\infty}(1+\xi^2)|x(\xi)|\,d\xi<\infty, \quad
\textup{b).}\ %
\int\limits_{0}^{\infty}\bigg|%
\frac{d\,}{d\xi}\bigg(\xi^2\frac{dx(\xi)}{d\xi}\bigg)\bigg|\,d\xi<\infty\,.
\end{equation}
We show, (see Lemma \ref{cfin} below), that under the conditions
\eqref{ssiic2},
 the terms out the integral in
\eqref{siip} vanish as \(\varepsilon\to{}0,\,N\to\infty\).
Therefore, under the conditions  \eqref{ssiic2},  we can pass to
the limit in \eqref{siip} as \(\varepsilon\to{}0,\,N\to\infty\):
\begin{equation}
\label{pts}
\int\limits_{0}^{\infty}\left(-\frac{d\,\,}{d\xi}\bigg(\xi^2
\frac{dx(\xi)}{d\xi}\bigg)\right)e^{it\xi}d\xi=%
-it\int\limits_{0}^{\infty}x(\xi)\Bigg(\frac{d\,}{d\xi}\,\bigg(\xi^2e^{it\xi}\bigg)\Bigg)\,d\xi\,.
\end{equation}
Transforming the integral in the right hand side of \eqref{pts}
and denoting
\begin{math}
 y(t)=\int\limits_{0}^{\infty}x(\xi)e^{it\xi}\,d(\xi),
\end{math}
 we obtain
 \begin{equation*}
 -it\int\limits_{0}^{\infty}x(\xi)\Bigg(\frac{d\,}{d\xi}\,\bigg(\xi^2e^{it\xi}\bigg)\Bigg)\,d\xi=
-t^2\frac{d^2y(t)}{dt^2}-2t\frac{dy(t)}{dt}=-\frac{d\,\,}{dt}\bigg(t^2\frac{dy(t)}{dt}\bigg)\,.
\end{equation*}
\begin{lemma}
\label{csfin} Let \(x(\xi)\) be a function defined on the open
semi-bounded interval \((0,\infty)\). We assume that:
\\
\hspace*{2.0ex}\textup{1.}
\begin{minipage}[t]{0.92\linewidth}
The function \(x(\xi)\) is differentiable on \((0,\infty)\) and
that its derivative \(\dfrac{dx(\xi)}{d\xi}\) is absolutely
continuous on every compact subinterval of the interval
\((0,\infty)\).
\end{minipage}\\[1.0ex]
\hspace*{2.0ex}\textup{2.}
\begin{minipage}[t]{0.92\linewidth}
The conditions \eqref{ssiic2} hold.
\end{minipage}\\[0.1ex]

Then the boundary conditions
\begin{subequations}
\label{BCSII}
\begin{alignat}{2}
\label{BCSIIa}%
 \lim_{\xi\to{}+0}\xi^2\,x(\xi)&=0,&\quad
\lim_{\xi\to{}+0}\xi^2\,\frac{dx(\xi)}{d\xi}&=0,
\\[1.0ex]
\label{BCSIIb}%
\lim_{\xi\to\infty}\xi^2\,x(\xi)&=0,&\quad
\lim_{\xi\to\infty}\xi^2\frac{dx(\xi)}{d\xi}&=0.
\end{alignat}
\end{subequations} are satisfied.
\begin{proof}
Let
\[0<\alpha<\beta<\infty\,.\]
Since
\[\xi^2\frac{dx(\xi)}{d\xi}\bigg|_{\xi=\beta}\!\!\!-
\xi^2\frac{dx(\xi)}{d\xi}\bigg|_{\xi=\alpha}\!\!\!=
\int\limits_{\alpha}^{\beta}%
\frac{d\,}{d\xi}\bigg(\xi^2 \frac{dx(\xi)}{d\xi}\bigg)\,d\xi,\] it
follows from (\ref{ssiic2}.b) that there exist finite limits
\[\lim_{\alpha\to+0}%
c_0=\xi^2\frac{dx(\xi)}{d\xi}\bigg|_{\xi=\alpha},\quad%
c_\infty=\lim_{\beta\to{}\infty}%
\xi^2\frac{dx(\xi)}{d\xi}\bigg|_{\xi=\beta}\,.\] Since
\(x(\xi_2)-x(\xi_1)=\int\limits_{\xi_1}^{\xi_2}\frac{dx(\xi)}{d\xi}d\xi\),
then
\[x(\xi)=-\frac{c_0}{\xi}+o(|\xi|^{-1})\ \ \textup{as}\ \ \xi\to+0\,.\]
Since \(\int\limits_0^1|x(\xi)|d\xi<\infty\), \(c_{{}_0}=0\).
Thus, the boundary conditions \eqref{BCSIIa} hold. In the same
way,
\[x(\xi)=-\frac{c_\infty}{\xi}+o(|\xi|^{-1})\ \ \textup{as}\ \ \xi\to\infty\,.\]
Since \(\displaystyle\int\limits_1^\infty|x(\xi)|d\xi<\infty\),
\(c_{\infty}=0\). Thus, the boundary condition %
\[\lim_{\xi\to{}\infty}%
\xi^2\frac{dx(\xi)}{d\xi}=0\] is satisfied.

The equality \(c_\infty=0\) ensure that
\(\displaystyle\lim_{\xi\to\infty}\bigg|\dfrac{d\,\,}{d\xi}(\xi^2x(\xi))\bigg|=0\,.\)
Thus, the function \(\xi^2x(\xi))\) satisfy the Lipschitz
condition at \(\infty\). (Even more, the Lipschitz constant equals
zero asymptotically.)) The convergence of the integral: %
\(\int\limits_0^\infty\xi^2|x(\xi)|d\xi<\infty\) together with the
Lipschitz condition for the function  \(\xi^2x(\xi))\), which
appears  under the integral, ensure that this function vanishes at
infinity. In other words,the boundary condition
\[\lim_{\xi\to{}\infty}%
\xi^2x(\xi)=0\] is satisfied.
\end{proof}
\end{lemma}
The above reasoning proves the following
\begin{theorem}
\label{csi}%
Let \(x(\xi)\) be a function defined on \((0,\infty)\).
 Assume that the following conditions are satisfied:\\
 \hspace*{2.0ex}\textup{1.}
\begin{minipage}[t]{0.92\linewidth}
The function \(x(\xi)\) is differentiable on \((0,a\infty)\) and
its derivative \(\dfrac{dx(\xi)}{d\xi}\) is absolutely continuous
on every compact subinterval of the interval \((0,\infty)\).
\end{minipage}\\[1.0ex]
\hspace*{2.0ex}\textup{2.}
\begin{minipage}[t]{0.92\linewidth}
\centerline{
\begin{math}
\displaystyle
\int\limits_{0}^{\infty}(1+\xi^2)|x(\xi)|\,d\xi<\infty, \quad
\int\limits_{0}^{\infty}\bigg|%
\frac{d\,}{d\xi}\bigg(\xi^2\frac{dx(\xi)}{d\xi}\bigg)\bigg|\,d\xi<\infty\,.
\end{math}}
\end{minipage}\\[1.0ex]

Then the function \(y(t)\):
\begin{equation}%
\label{FTiis}%
 y(t)=\int\limits_{0}^{\infty}x(\xi)\,e^{it\xi}\,d\xi,
\end{equation}
is defined for every \(t\in\mathbb{R}\), is twice continuously
differentiable on \(\mathbb{R}\), and the following commutation
relation holds:
\begin{equation}%
\label{CRsa}%
-\frac{d\,\,}{dt}\bigg(t^2\frac{dy(t)}{dt}\bigg)=
\int\limits_{0}^{\infty}\left(-\frac{d\,\,}{d\xi}\bigg(\xi^2
\frac{dx(\xi)}{d\xi}\bigg)\right)e^{it\xi}d\xi\,,\quad\forall\,
t\in\mathbb{R}\,.
\end{equation}
\end{theorem}
\begin{remark}
\label{RES}%
 In \textup{Theorems \ref{cwa}} and \textup{\ref{csi}}, no boundary conditions on
 the function \(x(\xi)\) were posed at the endpoints of the
interval. However, in \textup{Theorems \ref{cfi}} the boundary
conditions \eqref{bccfi} on the function \(x(\xi)\) were posed.
This is felt to be tied to the following fact. Let
\(\mathscr{L}_E^0\) be the \emph{minimal} differential operator in
the Hilbert space \(L^2(E)\) which is generated by the
differential expression \(L_E\) on the domain of definition
\(\mathcal{D}_0\), which consists of the set of all finite smooth
functions on \(E\). The minimal differential operator
\(\mathscr{L}_E^0\) is symmetric. In the cases
\(E=(-\infty,\infty)\) and \(E=(0,\infty)\), the operator
\(\mathscr{L}_E^0\) is an essentially selfadjoint operator on
\(\mathcal{D}_0\). This means that the closure
\(\overline{\mathscr{L}_E^0}\) of this operator is a selfadjoint
operator. However, in the case \(E={-a,a},\,0<a<\infty\), the
operator \(\overline{\mathscr{L}_E^0}\) is \emph{symmetric, but
not self-adjoint}. The deficiency indices of this operator are
equal to \((2,2)\). Selfadjoint extensions of this operator form
two-parametric family and may be described by boundary conditions
at the endpoints. Only one of these extensions commute with the
operator \(\mathscr{F}_E\). This is the extension that corresponds
to the boundary conditions \eqref{bccfi}. (We need to seek for
selfadjoint operators which commute with \(\mathscr{F}_E\) rather
then symmetric one. These are selfadjoint operators   which admit
the spectral decomposition.)
\end{remark}
\newpage
%%%%%%%%%%%%%%%%%%%%%%%%%%%%%%%%%%%%%%%%%%%%%%%%%%%%%%%%%%%%%%%%%%%%%%%%%%%%%%%%%%%%
%%%%%%%%%%%%%%%%%%%%%%%%%%%%%%%%%%%%%%%%%%%%%%%%%%%%%%%%%%%%%%%%%%%%%%%%%%%%%%%%%%%%%%%%

\vspace*{5.0ex}
\noindent
\begin{minipage}[h]{0.45\linewidth}
Victor Katsnelson\\[0.2ex]
Department of Mathematics\\
The Weizmann Institute\\
Rehovot, 76100, Israel\\[0.1ex]
e-mail:\\
{\small\texttt{victor.katsnelson@weizmann.ac.il}}
\end{minipage}
\vspace*{3.0ex}

\noindent
\begin{minipage}[h]{0.45\linewidth}
Ronny Machluf\\[0.2ex]
Department of Mathematics\\
The Weizmann Institute\\
Rehovot, 76100, Israel\\[0.1ex]
e-mail:\\
\texttt{ronny-haim.machluf@weizmann.ac.il}
\end{minipage}

\begin{thebibliography}{MoMo1}
%%%%%%%%%%%%%%%%%%%%%%%%%%%%%%%%%%%%%%%%%%%%%%%%%%%%%%%%%%%%%%%%%%%%%%%%%%%%%

%%%%%%%%%%%%%%%%%%%%%%%%%%%%%%%%%%%%%%%%%%%%%%%%%%%%%%%%%%%%%%%%%%%%%%%%%%%%
%%%%%%%%%%%%%%%%%%%%%%%%%%%%%%%%%%%%%%%%%%%%%%%%%%%%%%%%%%%%%%%%%%%%%%%%%%%
\bibitem[DuGr]{DuGr} \textsc{Duistermaat,\,J.J., Gr\"unbaum,\,F.A.}
\textit{Differential equations in spectral parameters.}
Communs.\,Math.\,Phys., \textbf{103} (1986), 177\,-\,240.
%%%%%%%%%%%%%%%%%%%%%%%%%%%%%%%%%%%%%%%%%%%%%%%%%%%%%%%%%%%%%%%%%%%%%%%%%%%%%%
\bibitem[GGK]{GGK}
\textsc{Gohberg,\,I.,\,Goldberg,\,S.,\,Kaashoek,M.A.}
\textit{Classes of Linear Operators. Vol.\,1.} Birkh\"auser,
Basel\(\boldsymbol{\cdot}\)Boston\(\boldsymbol{\cdot}\)Berlin
1990.
\bibitem[Gru1]{Gru1} \textsc{Gr\"unbaum, F.A.} \textit{Band-time-band limited
integral operators and commuting differential operators.} Àëãåáðà
è Àíàëèç, \textbf{8}:1 (1996), 122--126= St. Petersburg Math. J.,
\textbf{8}:1 (1997), 93\,-\,96.
%%%%%%%%%%%%%%%%%%%%%%%%%%%%%%%%%%%%%%%%%%%%%%%%%%%%%%%%%%%%%%%%%%%%%%%%%%%%%%
\bibitem[Gru2]{Gru2} \textsc{Gr\"unbaum, F.A.} \textit{Some bispectral musings.}
pp.\,31\,-\,45 in \cite{HaKa}.
%%%%%%%%%%%%%%%%%%%%%%%%%%%%%%%%%%%%%%%%%%%%%%%%%%%%%%%%%%%%%%%%%%%%%%%%%%%%%%%
\bibitem[GrYa]{GrYa}  \textsc{\mbox{Gr\"unbaum, F.A.,
Yakimov,\,M.}} \textit{The prolate spheroidal phenomenon as a
consequence of bispectrality.} pp.\,301\,-\,312 in \cite{TWHMPR}.
%%%%%%%%%%%%%%%%%%%%%%%%%%%%%%%%%%%%%%%%%%%%%%%%%%%%%%%%%%%%%%%%%%%%%%%%%%%%%
\bibitem[HaKa]{HaKa} \textsc{Harnad,\,J., Kasman,\,A.}-ed.
 \textit{The bispectral problem. {\footnotesize Papers from the CRM Workshop held at the
Université de Montreal, Montreal, PQ, March 1997}.} CRM
Proceedings \& Lecture Notes, 14. American Mathematical Society,
Providence, RI, 1998. viii+235 pp.
%%%%%%%%%%%%%%%%%%%%%%%%%%%%%%%%%%%%%%%%%%%%%%%%%%%%%%%%%%%%%%%%%%%%%%%%
\bibitem[LaP1]{LaP1} \textsc{Landau,\,H., Pollak,\,H.O.} \textit{Prolate spheroidal wave
functions, Fourier analysis and uncertainty\,--\,\textup{II}}.
Bell System Techn. Journ. \textbf{40} (1961), 65\,-\,84.
%%%%%%%%%%%%%%%%%%%%%%%%%%%%%%%%%%%%%%%%%%%%%%%%%%%%%%%%%%%%%%%%%%%%%%%%
\bibitem[LaP2]{LaP2} \textsc{Landau,\,H., Pollak,\,H.O.} \textit{Prolate spheroidal wave
functions, Fourier analysis and
uncertainty\,--\,\textup{III}:\,The dimension of the space of
essentially time- and band-limited signals}. Bell System Techn.
Journ. \textbf{40} (1961), 1295\,-\,1336.
%%%%%%%%%%%%%%%%%%%%%%%%%%%%%%%%%%%%%%%%%%%%%%%%%%%%%%%%%%%%%%%%%%%%%%%%
\bibitem[Mo1]{Mo1}\textsc{Morrison,\,J.A.} \textit{On the eigenfunctions
 corresponding to the bandpass kernel, in the case of degeneracy.}
 Quarterly of Appl. Math., \textbf{21} (1963), 13\,-\,19.
 %%%%%%%%%%%%%%%%%%%%%%%%%%%%%%%%%%%%%%%%%%%%%%%%%%%%%%%%%%%%%%%%%%%%
\bibitem[Mo2]{Mo2}\textsc{Morrison,\,J.A.} \textit{Eigenfunctions
of the finite Fourier transform operator over a hyperellipsoidal
region.} Journ. Math. and Mech., \textbf{44} (1965), 245\,-\,254.
%%%%%%%%%%%%%%%%%%%%%%%%%%%%%%%%%%%%%%%%%%%%%%%%%%%%%%%%%%%%%%%%%%%%%%%%%
\bibitem[Mo2]{Mo2}\textsc{Morrison,\,J.A.} \textit{Dual formulation for the
eigenfunctions correspondind to the bandpass kernel, in the case
of degeneracy.} Journ. Math. and Mech., \textbf{44} (1965),
313\,-\,326.
%%%%%%%%%%%%%%%%%%%%%%%%%%%%%%%%%%%%%%%%%%%%%%%%%%%%%%%%%%%%%%%%%%%%%%%%%%%%%%
\bibitem[Sl1]{Sl1} \textsc{Slepian,\,D.}
\textit{Prolated spheroidal wave functions, Fourier analysis and
uncertainity -- \textup{IV}: Extension to many dimensions;
generalized prolate spheroidal functions.} Bell System Techn.
Journ. \textbf{43} (1964), 3009\,-\,3057.
%%%%%%%%%%%%%%%%%%%%%%%%%%%%%%%%%%%%%%%%%%%%%%%%%%%%%%%%%%%%%%%%%%%%%%%%%%%%%%
\bibitem[Sl2]{Sl2} \textsc{Slepian,\,D.} \textit{Some comments on Fourier analysis, uncertainty
and modelling.} SIAM Review, \textbf{25}:3, 1983, 379\,-\,393.
%%%%%%%%%%%%%%%%%%%%%%%%%%%%%%%%%%%%%%%%%%%%%%%%%%%%%%%%%%%%%%%%%%%%%%%%%%%%%%
\bibitem[SlPo]{SlPo} \textsc{Slepian,\,D., Pollak,\,H.O.}
\textit{Prolated spheroidal wave functions, Fourier analysis and
uncertainity -- \textup{I}.} Bell System Techn. Journ. \textbf{40}
(1961), 43\,-\,63.
%%%%%%%%%%%%%%%%%%%%%%%%%%%%%%%%%%%%%%%%%%%%%%%%%%%%%%%%%%%%%%%%%%%%%%%%%%%%%%
\bibitem[TWHMPR]{TWHMPR} \textsc{\mbox{Tempesta,\,P., Winternitz,\,P., Harnad,\,J.,} %
Miller,\,W.,\,Jr., Pogosyan,\,G, Rodrigues,\,M.}\\
\textit{Superintegrability in Classical and Quantum Systems.} (CRM
Proc. \& Lect. Notes, Vol.\,37.) Amer. Math. Soc., Providence, RI,
2002, x+347 pp.
%%%%%%%%%%%%%%%%%%%%%%%%%%%%%%%%%%%%%%%%%%%%%%%%%%%%%%%%%%%%%%%%%%%%%%%%%%%%%%
\bibitem[Wal]{Wal} \textsc{Walter, G.} \textit{Differential operators which commute
with characteristic functions with application to a lucky
accident.} Complex Variables, \textbf{18} (1992), 7\,-\,12.
\end{thebibliography}
\end{document}